\documentclass[british]
{elsarticle}
\usepackage{amsmath, amssymb} 
\usepackage{babel,eucal,url,amssymb,
enumerate,amscd,
}

\begin{document}


\title
{A connection with parallel torsion 
on almost hypercomplex manifolds 
with Hermitian and
anti-Hermitian metrics\tnoteref{t2}} \tnotetext[t2]{Dedicated to
Professor Kostadin Gribachev, the author's teacher in the theory
of manifolds, on the occasion of his 70th birthday}

\author{Mancho Manev}
\address{Paisii Hilendarski University of Plovdiv \\
Faculty of Mathematics and Informatics\\
236 Bulgaria Blvd., 4003 Plovdiv, Bulgaria}
\ead{mmanev@uni-plovdiv.bg}



\newcommand{\ie}{i.~e. }
\newcommand{\X}{\mathfrak{X}}
\newcommand{\W}{\mathcal{W}}
\newcommand{\K}{\mathcal{K}}
\newcommand{\N}{\mathbb{N}}
\newcommand{\R}{\mathbb{R}}
\newcommand{\s}{\mathfrak{S}}
\newcommand{\n}{\nabla}
\newcommand{\D}{{\rm d}}
\newcommand{\al}{\alpha}
\newcommand{\bt}{\beta}
\newcommand{\gm}{\gamma}
\newcommand{\ea}{\varepsilon_\alpha}
\newcommand{\eb}{\varepsilon_\beta}
\newcommand{\eg}{\varepsilon_\gamma}
\newcommand{\sa}{\sum_{\al=1}^3}
\newcommand{\sbt}{\sum_{\bt=1}^3}
\newcommand{\ee}{\end{equation}}
\newcommand{\be}[1]{\begin{equation}\label{#1}}
\def\bea{\begin{eqnarray}} \def\eea{\end{eqnarray}}
\newcommand{\norm}[1]{\left\Vert#1\right\Vert ^2}
\newcommand{\nJ}[1]{\norm{\nabla J_{#1}}}

\newcommand{\thmref}[1]{The\-o\-rem~\ref{#1}}
\newcommand{\propref}[1]{Pro\-po\-si\-ti\-on~\ref{#1}}
\newcommand{\secref}[1]{\S\ref{#1}}
\newcommand{\lemref}[1]{Lem\-ma~\ref{#1}}
\newcommand{\dfnref}[1]{De\-fi\-ni\-ti\-on~\ref{#1}}
\newcommand{\coref}[1]{Corollary~\ref{#1}}


\newtheorem{thm}{Theorem}
\newtheorem{lem}[thm]{Lemma}
\newtheorem{prop}[thm]{Proposition}
\newtheorem{cor}[thm]{Corollary}
\newdefinition{rmk}{Remark}
\newdefinition{ack}{Acknowledgements}
\newproof{pf}{Proof}
\newproof{pot}{Proof of Theorem \ref{thm-geom}}




\hyphenation{Her-mi-ti-an ma-ni-fold ah-ler-ian}

\begin{keyword}
almost hypercomplex manifold \sep pseudo-Riemannian metric
 \sep anti-Hermitian metric \sep Norden metric \sep indefinite metric \sep neutral metric \sep
natural connection \sep parallel structure. \MSC{53C26, 53C15,
53C50, 53C55.}
\end{keyword}


\begin{abstract}
Almost hypercomplex manifolds with Hermitian and anti-Hermitian
metrics are considered. A linear connection $D$ is introduced such
that the structure of these manifolds is parallel with respect to
$D$. Of special interest is the class of the locally conformally
equivalent manifolds of the manifolds with covariantly constant
almost complex structures and the case when the torsion of $D$ is
$D$-parallel. Curvature properties of these manifolds are studied.
An example of 4-dimensional manifolds in the considered basic
class is constructed and characterized.
\end{abstract}

\maketitle



\section*{Introduction}

In this work\footnote{This work was partially supported by the
Scientific Researches Fund at the University of Plovdiv.} we
continue the investigations on
a manifold $M$ with an almost hypercomplex structure $H$
. We provide this almost hypercomplex manifold $(M,H)$ with a
metric structure $G$, generated by a pseudo-Riemannian metric
$g$ of neutral signature (\cite{GrMa}, \cite{GrMaDi}).

It is known that, if $g$ is a Hermitian metric on $(M,H)$, the
derived metric structure $G$ is the known hyper-Hermitian
structure. It consists of the given Hermitian metric $g$ with
respect to the three almost complex structures of $H$ and the
three K\"ahler forms associated with $g$ by $H$ \cite{AlMa}.
%

Here, the considered metric structure $G$ has a different type of
compatibility with $H$. The structure $G$ is generated by a
neutral metric $g$ such that the first (resp., the other two) of
the almost complex structures of $H$ acts as an iso\-metry (resp.,
act as anti-isometries) with respect to $g$ in each tangent
fibre. 
%
Let the almost complex structures of $H$ act as iso\-metries or
anti-isometries with respect to the metric, then the existence of
an anti-isometry generates exactly the existence of one more
anti-isometry and an isometry. %
Thus,
$G$ contains the metric $g$ and three (0,2)-tensors associated by
$H$ -- a K\"ahler form and two metrics
of the same type.
The existence of such bilinear forms on an almost hypercomplex
manifold is proved in \cite{GrMa}.
The neutral metric $g$ is Her\-mit\-ian with respect to the first
almost complex structure of $H$ and $g$ is an anti-Hermitian (\ie
 Norden) metric regarding the other two almost complex structures
of $H$. For this reason we call the derived manifold $(M,H,G)$ an
\emph{almost hypercomplex manifold with Hermitian and
anti-Hermitian metrics} or shortly an \emph{almost
$(H,G)$-manifold}.
%
%
%
%

Recently, manifolds with neutral metrics and various tensor
structures have been object of interest in theoretical physics.

The geometry of an arbitrary almost $(H,G)$-manifold is the
geometry of the hypercomplex structure $H=\{J_1,J_2,J_3\}$ and the
neutral metric $g$ or equivalently -- the geometry of the metric
structure $G=\{g,g_1,g_2,g_3\}$. In this geometry, there are
important so-called natural connections for the $(H,G)$-structure
(briefly, the $(H,G)$-connections), \ie those linear connections,
with respect to which the parallel transport determines an
isomorphism between the tangent spaces with $(H,G)$-structure.
This holds if and only if $H$ and $g$ are parallel with respect to
such a connection.

If the three almost complex structures of $H$ are parallel with
respect to the Levi-Civita connection $\nabla$ of $g$, then we
call such $(H,G)$-manifolds of K\"ahler type
\emph{pseudo-hyper-K\"ahler manifolds} and we denote their class
by $\K$. Therefore, outside of the class $\K$, the Levi-Civita
connection $\n$ is no longer an $(H,G)$-connection. There exist
countless natural connections on an almost $(H,G)$-manifold in the
general case.

The class $\W$ of the locally conformally equivalent manifolds  of
the pseudo-hyper-K\"ahler manifolds is an object of special
interest in this paper. The covariant derivatives of the elements
of $H$ are explicitly expressed  by the structure tensors of $H$
and $G$ on the manifolds in $\W$.

In this work we construct and characterize an $(H,G)$-connection
on $\W$-manifolds. In the case of almost hyper-Hermitian manifolds
such a connection is known as the Lichnerowicz connection
\cite{AlMa}.

In the first section we show that as the pseudo-hyper-K\"ahler
manifolds have a zero curvature tensor for $\nabla$
(\cite{GrMaDi}), thus the tensors with the same properties are
zero on any almost $(H,G)$-manifold.

In the second section we introduce a linear connection $D$ with
respect to which the structure tensors of the almost
$(H,G)$-manifolds are parallel.
Then we characterize the torsion tensor and the curvature tensor
of $D$.

In the third section we consider the special case when the torsion
of $D$ is $D$-parallel. There we characterize geometrically the
manifolds with that specialization.

In the fourth section we construct a class of 4-dimensional Lie
groups as $\W$-manifolds where the torsion of $D$ is not
$D$-parallel.

The basic problem of this work is the existence and the geometric
characteristics of the considered manifolds with $D$-parallel
torsion of a natural connection $D$. The main result of this paper
is that every locally conformal pseudo-hyper-K\"ahler manifold
with $D$-parallel torsion of $D$ is a $D$-flat Lie group.

\section{Almost hypercomplex manifolds with Hermitian and anti-Hermit\-ian metrics}

\subsection{The almost $(H,G)$-manifolds}

Let $(M,H)$ be an almost hypercomplex manifold, \ie $M$ is a
$4n$-dimension\-al differentiable manifold and $H=(J_1,J_2,J_3)$
is a triple of almost complex structures with the properties:
\be{J123} J_\al=J_\bt\circ J_\gm=-J_\gm\circ J_\bt, \qquad
J_\al^2=-I\ee for all cyclic permutations $(\al, \bt, \gm)$ of
$(1,2,3)$ and $I$ denotes the identity.

The standard structure of $H$ on a $4n$-dimensional vector space
with a basis
$\{X_{4k+1},X_{4k+2},X_{4k+3},X_{4k+4}\}_{k=0,1,\dots,n-1}$ has
the form \cite{So}:
\begin{equation}\label{Jdim4n}
\begin{array}{lll}
J_1X_{4k+1}= X_{4k+2},\quad & J_2X_{4k+1}=X_{4k+3},\quad & J_3X_{4k+1}=-X_{4k+4},\\[4pt]
J_1X_{4k+2}=-X_{4k+1},\quad & J_2X_{4k+2}=X_{4k+4}, \quad & J_3X_{4k+2}=X_{4k+3}, \\[4pt]
J_1X_{4k+3}=-X_{4k+4},\quad & J_2X_{4k+3}=-X_{4k+1}, \quad & J_3X_{4k+3}=-X_{4k+2},\\[4pt]
J_1X_{4k+4}=X_{4k+3},\quad & J_2X_{4k+4}=-X_{4k+2}, \quad &
J_3X_{4k+4}=X_{4k+1}.\\[4pt]
\end{array}
\end{equation}


Let $g$ be a pseudo-Riemannian metric on $(M,H)$ with the
properties
\be{gJJ} %
g(x,y)=\ea g(J_\al x,J_\al y), \ee %
where %
\be{ea}%
 \ea=
\begin{cases}
\begin{array}{ll}
1, \quad & \al=1;\\[4pt]
-1, \quad & \al=2;3.
\end{array}
\end{cases}
\ee
In other words, for $\al=1$, the metric $g$ is Hermitian with
respect to $J_1$, and in the case  $\al=2$ or $\al=3$, the metric
$g$ is an anti-Hermitian (\ie Norden) metric with respect to
$J_\al$ ($\al=2$ or $\al=3$, respectively) \cite{GaBo}. Moreover,
the associated bilinear forms $g_1$, $g_2$, $g_3$ are determined
by
\be{gJ} g_\al(x,y)=g(J_\al x,y)=-\ea g(x,J_\al y),\qquad
\al=1,2,3. \ee
Because of \eqref{gJJ} and \eqref{gJ}, the metric $g$ and the
associated bilinear forms $g_2$ and $g_3$ are necessarily
pseudo-Riemannian metrics of neutral signature $(2n,2n)$. The
associated bilinear form $g_1$ is the associated (K\"ahler)
2-form.

A structure $(H,G)=(J_1,J_2,J_3,g,g_1,g_2,g_3)$ is introduced in
\cite{GrMaDi} and \cite{GrMa}. The cases when the original metric
$g$ is Hermitian or anti-Hermitian with respect to the almost
complex structures of $H$ are considered. In \cite{GrMa} it is
proved that the unique possibility that an anti-Hermitian metric
be considered on an almost hypercomplex manifold is the case when
the given metric is Hermitian with respect to the first and
moreover it is an anti-Hermitian metric with respect to other two
structures of $H$. Therefore, we call $(H,G)$ an \emph{almost
hypercomplex structure with Hermitian and anti-Hermitian metrics}
on $M$ (or, in short, an \emph{almost $(H,G)$-structure} on $M$).
Then, we call briefly a manifold with such a structure an
\emph{almost $(H,G)$-manifold}.

The structural tensors of an almost $(H,G)$-manifold are the three
$(0,3)$-tensors determined by
\begin{equation}\label{F}
F_\al (x,y,z)=g\bigl( \left( \n_x J_\al
\right)y,z\bigr)=\bigl(\n_x g_\al\bigr) \left( y,z \right),\qquad
\al=1,2,3,
\end{equation}
where $\n$ is the Levi-Civita connection generated by $g$.
The corresponding Lee forms $\theta_\al$ are defined by
\begin{equation}\label{theta-al}
\theta_\al(\cdot)=g^{ij}F_\al(e_i,e_j,\cdot),\qquad \al=1,2,3,
\end{equation}%
for an arbitrary basis $\{e_i\}_{i=1}^{4n}$.

The tensors $F_\al$ have the following fundamental identities:
\be{34}%
\begin{array}{c}
    F_\al(x,y,z)=-\ea F_\al(x,z,y)=-\ea F_\al(x,J_\al y,J_\al
    z),\\[6pt]
    F_\al(x,J_\al y,z)=\ea F_\al(x,y,J_\al z);
\end{array}
\ee
%
\begin{equation}\label{F-prop-1}
\begin{array}{l}
    F_\al(x,y,z)=F_\bt(x,J_\gm y,z)-\eb F_\gm(x,y,J_\bt z)\\[6pt]
    \phantom{F_\al(x,y,z)}=-F_\gm(x,J_\bt y,z)+\eg F_\bt(x,y,J_\gm
    z);
\end{array}
\end{equation}
\begin{equation}\label{F-prop-2}
\begin{array}{l}
    F_\bt(x,J_\gm y,z)-\eg F_\bt(x,y,J_\gm z)\\[6pt]
    \phantom{F_\bt(x,J_\gm y,z)}+F_\gm(x,J_\bt y,z)-\eb F_\gm(x,y,J_\bt
    z)=0;
\end{array}
\end{equation}
\be{F1ab}%
\begin{array}{l} %
F_\al(x,J_\bt y, J_\gm z)=\ea F_\al(x,J_\gm y,J_\bt z), \\[6pt]%
F_\al(x,J_\bt y, J_\bt z)=-\ea F_\al(x,J_\gm y,J_\gm z)%
\end{array} %
\ee %
for all cyclic permutations $(\al, \bt, \gm)$ of $(1,2,3)$.

In \cite{GrMaDi} we study a special class $\K$ of the
$(H,G)$-mani\-folds -- the so-called there
\emph{pseudo-hyper-K\"ahler manifolds}. The manifolds from the
class $\K$ are the $(H,G)$-mani\-folds for which the complex
structures $J_\al$ are parallel with respect to the Levi-Civita
connection $\n$, generated by $g$, for all $\al=1,2,3$.
%

As $g$ is an indefinite metric, there exist isotropic vectors $x$
on $M$, \ie \(g(x,x)=0\) for a nonzero vector \(x\). 
In \cite{GrMa} we define the invariant square norm
\begin{equation}\label{nJ}
\nJ{\alpha}= g^{ij}g^{kl}g\bigl( \left( \nabla_i J_\alpha \right)
e_k, \left( \nabla_j J_\alpha \right) e_l \bigr),
\end{equation}
where $\{e_i\}_{i=1}^{4n}$ is an arbitrary basis of the tangent
space $T_pM$ at an arbitrary point $p\in M$. We say that an almost
$(H,G)$-manifold is an \emph{isotropic pseudo-hyper-K\"ahler
manifold} if $\nJ{\alpha}=0$ for all $J_\alpha$ of $H$. Clearly,
if $(M,H,G)$ is a pseudo-hyper-K\"ahler manifold, then it is an
isotropic pseudo-hyper-K\"ahler manifold. The inverse statement
does not hold. For instance, in \cite{GrMa} we have constructed an
almost $(H,G)$-manifold on a Lie group, which is an isotropic
pseudo-hyper-K\"ahler manifold but it is not a
pseudo-hyper-K\"ahler manifold.

\subsection{Properties of the K\"ahler-like tensors}

A tensor $L$ of type (0,4) with the pro\-per\-ties:%
\be{curv}%
\begin{split}%
&L(x,y,z,w)=-L(y,x,z,w)=-L(x,y,w,z),\\[6pt]
&L(x,y,z,w)+L(y,z,x,w)+L(z,x,y,w)=0
\end{split}%
\ee %
is called a \emph{curvature-like tensor}. The last equality of
\eqref{curv} is known as the first Bianchi identity of a
curvature-like tensor $L$.

We say that a curvature-like tensor $L$ is a \emph{K\"ahler-like
tensor} on an almost $(H,G)$-manifold when $L$ satisfies the properties: %
\be{L-kel}%
\begin{array}{l}
L(x,y,z,w)=\ea L(x,y,J_\al z,J_\al w)=\ea L(J_\al x,J_\al y,z,w),
\end{array}
\ee where $\ea$ is determined by \eqref{ea}.

Let the curvature tensor $R$ of the Levi-Civita connection
$\nabla$, generated by $g$, be defined, as usual, by
$R(x,y)z=\nabla_x \nabla_y z - \nabla_y \nabla_x z -
\nabla_{\left[x,y\right]} z$. The corresponding $(0,4)$-tensor,
denoted by the same letter, is determined by
$R(x,y,z,w)=g\left(R(x,y)z,w\right)$. Obviously, $R$ is a
K\"ahler-like tensor on an arbitrary pseudo-hyper-K\"ahler
manifold.

A K\"ahler-like tensor $L$ on an arbitrary almost $(H,G)$-manifold
has the same properties \eqref{curv} and \eqref{L-kel} of $R$ on a
pseudo-hyper-K\"ahler manifold. Then, according to \cite{GrMaDi},
we have the following similar properties:
\[
L(x,y,z,w)=\ea L(x,J_\al y,J_\al z,w),
\]
\[
L(x,y,z,w)=-L(x,J_\al y,z,J_\al w)=-L(J_\al x,y,J_\al z,w).
\]

Thus, we obtain the following geometric characteristic of the
K\"ahler-like tensors on an almost $(H,G)$-manifold, similarly to
Theorem~2.3 in \cite{GrMaDi}, where it is proved that the
hyper-K\"ahler $(H,G)$-manifolds are flat, \ie $R=0$ in $\K$.
\begin{prop}\label{th-0}
Every K\"ahler-like tensor on an almost $(H,G)$-manifold is
ze\-ro.  \phantom{} \hfill $\Box$
\end{prop}

For the sake of brevity we introduce the following notation for
an arbitrary (0,4)-tensor  %
\[%
(L\circ J_\al)(x,y,z,w)=L(x,y,J_\al
z,J_\al w).%
\]

Let us recall the Ricci identity for an almost complex structure
$J$ and a curvature tensor $R$ for the Levi-Civita connection $\n$
of $g$:
\[
\left(\n_x\n_y J\right)z-\left(\n_y\n_x
J\right)z=R(x,y)Jz-JR(x,y)z.
\]
Then, according to \eqref{F}, \eqref{gJ}, \eqref{ea}, \eqref{gJJ}
and $\n g=0$, we obtain the following corollaries of the Ricci
identity for the almost complex structures $J_\al$ $(\al=1,2,3)$ of $H$:%
\be{Ric-id} %
\begin{split}%
&\left(\n_x F_\al\right)(y,z,w)-\left(\n_y
F_\al\right)(x,z,w)\\[6pt]
&\phantom{\left(\n_x F_\al\right)(y,z,w)}=R(x,y,J_\al z,w)+\ea
R(x,y,z,J_\al w),
\end{split}
\ee
\be{Ric-id-J} %
\begin{split}%
&R(x,y,z,w)-\ea R(x,y,J_\al z,J_\al w)\\[6pt]
&\phantom{R(x,y,z,w)} =-\ea\bigl\{\left(\n_x
F_\al\right)(y,z,J_\al w)-\left(\n_y F_\al\right)(x,z,J_\al
w)\bigr\} \\[6pt]%
&\phantom{R(x,y,z,w)} =-\left(\n_x F_\al\right)(y,J_\al z,
w)+\left(\n_y F_\al\right)(x,J_\al z,w).
\end{split}
\ee

\subsection{The class $\W$ of the locally conformal pseudo-hyper-K\"ahler manifolds}

According to \eqref{gJJ} for $\alpha=1$, the manifold $(M,J_1,g)$
is almost Hermitian. The following class $\W(J_1)$ is the class
where the tensor $F_{1}$ is expressed explicitly by the structure
tensors. 
It is
denoted by $\W_4$ in \cite{GrHe}. For dimension $4n$ this class is
determined by
\begin{equation}\label{cl-H}
\begin{split}
        \W(J_1):\; F_1(x,y,z)=&
        \frac{1}{2(2n-1)}
                \left\{g(x,y)\theta_1(z)-g(x,z)\theta_1(y)\right.
                \\[6pt]
                &
                \left.
                -g(x,J_1y)\theta_1(J_1z)+g(x,J_1z)\theta_1(J_1y)
                \right\}.
\end{split}
\end{equation}

On the other hand, having in mind \eqref{gJJ} for $\alpha=2$ or
$3$, the manifold $(M,J_\alpha,g)$ is an almost complex manifold
with Norden metric (\ie anti-Hermitian or $B$-metric). The class
$\W_1$ in \cite{GaBo} is the class where the tensor $F_{\al}$ is
expressed explicitly by the structure tensors. We denote this
class by $\W(J_\alpha)$. For dimension $4n$ and \(\al =2,3\), it
is determined by
\begin{equation}\label{cl-N}
\begin{split}
\W(J_\al):\; F_\al(x,y,z)=& \frac{1}{4n}\bigl\{
g(x,y)\theta_\al(z)+g(x,z)\theta_\al(y)\bigr.\\[6pt]
& \bigl.+g(x,J_\al y)\theta_\al(J_\al z)
    +g(x,J_\al z)\theta_\al(J_\al y)\bigr\}.\\[6pt]
\end{split}
\end{equation}

The three special classes $\W_0(J_\al)$: $F_\al=0$ for $\al=1,2,3$
of the K\"ahler-type manifolds belong to any other class within
the corresponding classification.

Let us denote the class $\W=\bigcap_{\alpha=1}^3\W(J_\alpha)$ of
the considered manifolds. It is known from \cite{GrMaDi}, that if
a manifold $(M,H,G)$ belongs to the class  $ \W(J_\alpha) \bigcap
\W(J_\beta)$, then $(M,H,G)$ is of the class $\W(J_\gamma)$ for
all cyclic permutations $(\alpha, \beta, \gamma)$ of $(1,2,3)$.

It is well known that the almost hypercomplex structure
$H=(J_\al)$ is a {\em hypercomplex structure\/} if $N_\al$
vanishes for all $\al=1,2,3$, where \( N_\al(\cdot,\cdot)=
\left[J_\al \cdot,J_\al \cdot \right]
    -J_\al\left[J_\al \cdot,\cdot \right]
    -J_\al\left[\cdot,J_\al \cdot \right]
    -\left[\cdot,\cdot \right]\) are the Nijenhuis tensors for $J_\al$.
Moreover, it is known that an almost hypercomplex structure $H$ is
hypercomplex if and only if two of the three structures $J_\al$
are integrable. This means that two of the three tensors $N_\al$
vanish \cite{AlMa}. According to \cite{GrMaDi}, the class $\W$ is
a subclass of the class of the (integrable) $(H,G)$-manifolds. An
$(H,G)$-manifold belonging to $\W$ will be called in brief a
\emph{$\W$-manifold}.

Let us recall from \cite{GrMaDi} that if $(M,H,G)$ belongs to the
class $\W$ then
\[
\frac{2n}{1-2n}\theta_1\circ J_1=\theta_2\circ J_2=\theta_3\circ
J_3,
\]
where $\theta_\al$ $(\al=1,2,3)$ are defined in \eqref{theta-al}.
Hence, having in mind the last equalities, we introduce an 1-form
$\theta$ as follows
\begin{equation}\label{titi-all}
\theta=\frac{4n}{1-\ea(4n-1)}\theta_\al\circ J_\al,
\end{equation}
where $\al=1,2,3$ and $\ea$, $\theta_\al$ are determined by
\eqref{ea}, \eqref{theta-al}, respectively.

Let $c: g \mapsto \bar{g} = e^{2u}g$, where $u$ is a differential
function on $M$, be the usual conformal transformation of the
given metric. Then the manifolds  $(M,H,G)$ and $(M,H,\bar{G})$,
generated by the metrics $g$ and $\bar{g}$, respectively, are
called \emph{locally conformally equivalent manifolds}. Let
$(M,H,G)$ be an arbitrary pseudo-hyper-K\"ahler manifold. Then we
consider the class of the locally conformally equivalent manifolds
$(M,H,\bar{G})$ to the pseudo-hyper-K\"ahler manifolds (\ie in
short, locally conformal pseudo-hyper-K\"ahler manifolds).

\begin{thm}\label{th-conf}
The class $\W$ is the class of the locally conformally equivalent
manifolds to the pseudo-hyper-K\"ahler manifolds.
\end{thm}
\begin{pf}
Let $c: \bar g = e^{2u}g$ be the usual conformal transformation,
where $u$ is a differential function on $M$. Then, according to
Theorem 3.2 in \cite{GrMaDi}, the subclass $\W^0$ of $\W$ with the
condition $\D\theta=0$ is the class of the locally conformally
equivalent manifolds of the pseudo-hyper-K\"ahler manifolds.

According to \cite{GaGrMi}, the 1-form $\theta_\al\circ J_\al$ is
closed on the complex manifold with anti-Hermitian metric
$(M,J_\al,g)$ in $\W(J_\al)$ ($\al=2,3$) for any even dimension of
$M$. Therefore, having in mind \eqref{titi-all}, the 1-forms
$\theta_1\circ J_1$ and $\theta$ are closed for any
$4n$-dimensional $\W$-manifold ($n\in\N$). Hence, $\W^0$ coincides
with $\W$. \hfill$\Box$\end{pf}

Let us remark that for a Hermitian manifold $(M,J_1,g)$ in
$\W(J_1)$ the 1-form $\theta_1\circ J_1$ is closed for a dimension
greater than 4.

From the last theorem and the fact that every
pseudo-hyper-K\"ahler manifold is flat (\cite{GrMaDi}), we have
\begin{cor}\label{cor-conf-flat}
Every locally conformal pseudo-hyper-K\"ahler manifold is
conformally flat.$\hfill\Box$
\end{cor}

By virtue of \eqref{cl-H}, \eqref{cl-N} and \eqref{titi-all}, the
tensors $F_\al$ have the following form on a $\W$-manifold:
\begin{equation}\label{F123}
  \begin{array}{ll}
    &F_\al(x,y,z)=\frac{1}{4n}\left\{\ea g(x,y)\theta(J_\al z)-g(x,z)\theta(J_\al y)\right.
    \\[6pt]
                &\phantom{F_1(x,y,z)=\frac{1}{4n}.}
    \left.-\ea g(J_\al x,y)\theta(z)+g(J_\al
    x,z)\theta(y)\right\}.
  \end{array}
\end{equation}

On a $\W$-manifold the relations in \eqref{titi-all} between the
1-forms $\theta$ and $\theta_\al$ $(\al=1,2,3)$ imply immediately
corresponding equalities for the squares of the corresponding
vectors with respect to $g$ given in the following
\begin{cor}\label{cor-tOm}
On a locally conformal pseudo-hyper-K\"ahler manifold the
following equalities are valid
\[
\theta(\Omega)=\frac{\ea
16n^2}{[1-\ea(4n-1)]^2}\theta_\al(\Omega_\al),
\]
where $\Omega$ and $\Omega_\al$ are the Lee vectors associated to
$\theta$ and $\theta_\al$, respectively, \ie
$\theta=g(\cdot,\Omega)$ and $\theta_\al=g(\cdot,\Omega_\al)$.
Obviously, if one of the Lee vectors is isotropic then all Lee
vectors and their $J_\al$-images are isotropic.$\hfill\Box$
\end{cor}

According to \eqref{nJ}, \eqref{F}, \eqref{F123} and
\coref{cor-tOm}, we calculate immediately the square norms of
$\nabla J_\al$ on a $\W$-manifold.
\begin{prop}\label{prop-nJ}
On a locally conformal pseudo-hyper-K\"ahler manifold the squ\-are
norms of $\n J_\al$ for all $\al=1,2,3$ are proportional to the
square of the Lee vector $\Omega$, corresponding to the Lee form
$\theta$ from \eqref{titi-all}, with respect to $g$. More
precisely,
\[
\nJ{\al}=\frac{(4n-1)\ea-1}{4n^2}\theta(\Omega),\qquad \al=1,2,3.
\]
\hfill $\Box$
\end{prop}
Therefore we have
\begin{cor}\label{cor-iK}
A locally conformal pseudo-hyper-K\"ahler manifold is an isotropic
pseudo-hyper-K\"ahler manifold if and only if the Lee vector
$\Omega$ is isotropic.$\hfill \Box$
\end{cor}


\section{A natural $(H,G)$-connection}

\subsection{Definition of the connection $D$}
Let $(M,H,G)$ be an almost $(H,G)$-manifold generated by
$H=(J_1,J_2,J_3)$ and the metric $g$. We consider a linear
con\-nection on $(M,H,G)$ determined by
\be{D}
    D_y z=\n_y z+Q(y,z),\qquad Q(y,z)=\frac{1}{4}\sa \left(\n_y J_\al\right)J_\al
    z,
\ee %
 where $\n$ is the Levi-Civita connection generated by $g$. %
Obviously, the connections $D$ and $\n$ coincide if and only if
the manifold belongs to the class $\K$ of the
pseudo-hyper-K\"ahler manifolds.

By direct computations we establish that $D$ is a natural
connection for the structure $(H,G)$ or the following is valid
\begin{prop}\label{prop-natur}
The connection $D$ preserves the almost $(H,G)$-structure, \ie
\be{DH} D J_\al=Dg=Dg_\al=0,\qquad \al=1,2,3.\ee $\hfill \Box$
\end{prop}


Having in mind \propref{prop-natur}, we call the natural
connection $D$, determined by \eqref{D}, an
\emph{$(H,G)$-connection}.

Further, we use the corresponding tensor $Q$ of type (0,3)
determined by
\be{Q} %
Q(y,z,w)=g\left(Q(y,z),w\right)=\frac{1}{4}\sa
F_\al(y,J_\al z,w). %
\ee
Hence, according to \eqref{34}, we obtain %
\be{Q-Q}%
Q(x,y,z)=-Q(x,z,y). %
\ee

Since $DJ_\al=0$, we have $\left(\n_yJ_\al\right)z=-Q(y,J_\al
z)+J_\al Q(y,z)$ and consequently
\be{FQ}%
 F_\al(x,y,z)=-Q(x,J_\al y,z)-\ea Q(x,y,J_\al z), \qquad \al=1,2,3.
\ee
Therefore, because of \eqref{Q}, the following property of $Q$ is
valid:
\be{Q-prop}%
Q(x,y,z)+\sa \ea Q(x,J_\al y,J_\al z)=0.%
\ee

\subsection{The torsion tensor of $D$}

For the torsion $T$ of an arbitrary linear connection $D=\n+Q$, where $\n$ is the Levi-Civita connection, we have%
\be{T}
    T(x,y,z)=Q(x,y,z)-Q(y,x,z).
\ee
By virtue of \eqref{Q-Q} and \eqref{T} we obtain %
\be{QT}%
Q(x,y,z)=\frac{1}{2}\left\{T(x,y,z)-T(y,z,x)+T(z,x,y)\right\}. %
\ee

If the considered $(H,G)$-manifold belongs to the class $\W$,
using \eqref{T}, \eqref{Q} and \eqref{F123}, we obtain the
following
explicit form of the torsion of $D$ on a $\W$-manifold: %
\be{T=}
\begin{split}
    &T(x,y,z)=\frac{1}{16n}\bigl\{-2g_1(x,y)\theta(J_1z)\bigr.
+3g(x,z)\theta(y)-3g(y,z)\theta(x)\\[6pt]
&\phantom{T(x,y,z)=\frac{1}{16n}\bigl\{\bigr.}
+g_1(x,z)\theta(J_1y)+g_2(x,z)\theta(J_2y)+g_3(x,z)\theta(J_3y)\\[6pt]
&\phantom{T(x,y,z)=\frac{1}{16n}\bigl\{\bigr.}
-g_1(y,z)\theta(J_1x)-g_2(y,z)\theta(J_2x)-g_3(y,z)\theta(J_3x)\bigr\}.
\end{split}
\ee %
Hence we establish directly that the torsion of $D$ on a
$\W$-manifold has the property
\be{T-om=0} %
T(\cdot,\cdot,\Omega)=0,
\ee %
where $\Omega$ is the Lee vector, corresponding to the Lee form
$\theta$ from \eqref{titi-all}, with respect to $g$.

\subsection{The curvature tensor of $D$}
Let us consider the curvature tensor $K$ of a connection $D$, \ie
$K(x,y)z=[D_x, D_y] z - D_{\left[x,y\right]} z$ and
$K(x,y,z,w)=g\left(K(x,y)z,w\right)$. The relation between the
connections $D$ and $\n$ generates the corresponding relation
between their curvature tensors $K$ and $R$.
\begin{prop}\label{prop-KR}
On an almost $(H,G)$-manifold the curvature tensors $K$ of the
$(H,G)$-connection $D$, determined by \eqref{D}, and $R$ of
the Levi-Civita connection $\n$ are related as follows %
\be{K} %
K=\frac{1}{4}\{ R+R\circ J_1-R\circ J_2-R\circ
J_3\}-\frac{1}{16}P, %
\ee %
where %
\be{P} %
P=\sa\sbt P_{\al \bt}-4\sa P_{\al \al} %
\ee %
and
\[
\begin{array}{ll}
P_{\al \bt}(x,y,z,w)=&g\bigl(\left(\n_xJ_\al\right)J_\al z,
\left(\n_yJ_\bt\right)J_\bt
w\bigr)\\[6pt]
&-g\bigl(\left(\n_yJ_\al\right)J_\al z,
\left(\n_xJ_\bt\right)J_\bt w\bigr).
\end{array}%
\]
\end{prop}

\begin{pf}
For an arbitrary linear connection $D$ given by $D=\n+Q$ the
following equality is known \cite{Ko-No}%
\be{KR}
\begin{split}
&K(x,y,z,w)=R(x,y,z,w)+\left(\n_x Q\right)(y,z,w)-\left(\n_y
Q\right)(x,z,w)
\\[6pt]
&\phantom{K(x,y,z,w)=R(x,y,z,w)}+Q\bigl(x,Q(y,z),w\bigr)-Q\bigl(y,Q(x,z),w\bigr).
\end{split}
\ee %
Applying \eqref{Q} for the considered connection $D$, we obtain
consecutively  %
\be{QQ}
\begin{array}{l}
Q\bigl(x,Q(y,z),w\bigr)-Q\bigl(y,Q(x,z),w\bigr)=\\[6pt]
= g\bigl(Q(x,z),
Q(y,w)\bigr)-g\bigl(Q(y,z),Q(x,w)\bigr)\\[6pt]
=\frac{1}{16}\sa \sbt P_{\al \bt}(x,y,z,w),
\end{array}
\ee
\be{nQ1}%
\begin{split}
&\left(\n_x Q\right)(y,z,w)=\frac{1}{4}\sa \bigl\{\left(\n_x
F_\al\right)(y,J_\al z,w)\bigr.\\
&\phantom{\left(\n_x Q\right)(y,z,w)=\frac{1}{4}\sa\bigl\{\bigr.}
\bigl.-\ea g\bigl( \left(\n_x J_\al\right)z,\left(\n_y
J_\al\right)w\bigr)\bigr\}
\end{split}
\ee and therefore by \eqref{Ric-id-J} we have
\be{nQ}
\begin{split}
&\left(\n_x Q\right)(y,z,w)-\left(\n_y Q\right)(x,z,w)=
\\[6pt]
&\phantom{\left(\n_x Q\right)(y,z,w)}=\frac{1}{4}\sa \left\{
[-R+\ea (R\circ J_\al)-P_{\al \al}](x,y,z,w)\right\}.
\end{split}
\ee
We replace \eqref{nQ} and \eqref{QQ} in \eqref{KR} and then we
obtain \eqref{K}. \hfill$\Box$\end{pf}

When we consider $\W$-manifolds instead of arbitrary almost
$(H,G)$-manifolds, we have to specify the tensor $P$ introduced in
\eqref{P}.

Further we use the notation $g\odot h$ for the Kulkarni-Nomizu
product of two (0,2)-tensors, \ie
\[\begin{array}{ll} &\left(g\odot
h\right)(x,y,z,w)=g(x,z)h(y,w)-g(y,z)h(x,w)\\[6pt]
&\phantom{\left(g\odot
h\right)(x,y,z,w)=}+g(y,w)h(x,z)-g(x,w)h(y,z).\end{array}\]
Obviously, the tensor $g\odot h$ is a curvature-like tensor when
the (0,2)-tensors $g$ and $h$ are symmetric.

\begin{prop}\label{prop-W-P-2}
On a locally conformal pseudo-hyper-K\"ahler manifold the
curvature tensors $K$ and $R$ of the $(H,G)$-connection $D$ and
the Levi-Civita connection $\n$, respectively, are related by
\eqref{K}
and the auxiliary tensor $P$ has the form%
\be{PV} %
P=\frac{1}{16n^2}\{V+V\circ J_1-V\circ J_2-V\circ J_3\},%
\ee %
where %
\begin{equation}\label{barB}%
\begin{array}{ll}
&V=g\odot B+\frac{1}{2}U,\\[6pt]
& B=3\theta\otimes\theta+\theta\circ J_1\otimes\theta\circ
J_1+\theta\circ J_2\otimes\theta\circ J_2+\theta\circ
J_3\otimes\theta\circ
J_3\\[6pt]
&\phantom{
B=}-\frac{3}{2}\theta(\Omega)g-\frac{1}{2}\theta(J_2\Omega)g_2
-\frac{1}{2}\theta(J_3\Omega)g_3,\\[6pt]
&U(x,y,z,w)=g_1(x,y)E(z,w), \\[6pt]
&E=\theta\otimes\theta\circ J_1-\theta\circ
J_1\otimes\theta-\theta\circ J_2\otimes\theta\circ J_3+\theta\circ
J_3\otimes\theta\circ J_2.
\end{array}
\end{equation}
\end{prop}

\begin{pf}
The introduced above tensors have the following properties:
\begin{enumerate}\renewcommand{\labelenumi}{(\roman{enumi})}
    \item The (0,2)-tensor $ B$ is symmetric and $ B(J_\al \cdot,J_\al \cdot)\neq\pm  B(\cdot,\cdot)$;
    \item The (0,2)-tensor $E$ is antisymmetric and $E(J_\al \cdot,J_\al \cdot)=\ea
E(\cdot,\cdot)$;
    \item The (0,4)-tensor $U$ has the following properties: %
\be{U}
\begin{array}{c} %
U(x,y,z,w)=-U(y,x,z,w)=-U(x,y,w,z),\\[6pt]
U(x,y,J_\al z,J_\al w)=\ea U(x,y,z,w); %
\end{array}%
\ee %
    \item The (0,4)-tensor $V$ shares
the properties \eqref{U}.
\end{enumerate}

Equations \eqref{P} and \eqref{F123} imply directly %
\be{PW} P=\frac{1}{16n^2}\left\{(g\odot
B)+\sum_{\al=1}^3\ea(g\odot B)\circ J_\al\right\}+2U. %
\ee %

By virtue of the second line in \eqref{U}, we have %
\be{UU}%
U=\frac{1}{4}\{U+U\circ J_1-U\circ J_2-U\circ J_3 \}. %
\ee %

Then, from \eqref{PW} and \eqref{UU} we obtain \eqref{PV}. %
\hfill$\Box$\end{pf}

Now, let us consider the tensor %
\be{S}
S=R-\frac{1}{64n^2}\bigl(g\odot B\bigr), %
\ee %
where $B$ is determined by \eqref{barB}. Since $ B$ is symmetric,
then $S$ is a curvature-like tensor, where $S\circ J_\al\neq\pm
S$. Having in mind \propref{prop-KR} and \propref{prop-W-P-2}, we
have

\begin{prop}\label{prop-K-barS}
On a locally conformal pseudo-hyper-K\"ahler manifold the
curvature tensor
$K$ of the $(H,G)$-connection $D$ has the form%
\be{KSS} K=\frac{1}{4}\{\hat{S}+\hat{S}\circ J_1-\hat{S}\circ
J_2-\hat{S}\circ J_3\},\ee where
\[ \hat{S}=S-\frac{1}{128n^2}U,
\]
and $S$, $U$ are determined by \eqref{S}, \eqref{barB},
respectively.
\end{prop}

\begin{pf}
From \eqref{K} and \eqref{PV}, according to \eqref{S} and
\eqref{barB}, we have
\be{KS} %
K=\frac{1}{4}\{S+S\circ J_1-S\circ J_2-S\circ
J_3-\frac{1}{32n^2}U \}. %
\ee %
Then, taking into account \eqref{UU}, we obtain
 \eqref{KSS}. %
\hfill$\Box$\end{pf}

\section{Almost $(H,G)$-manifolds with $D$-parallel torsion of $D$}\label{D-parallel section}

\subsection{The condition that the torsion of $D$ be $D$-parallel}
In \secref{D-parallel section} we  consider the special case when
the torsion $T$ of $D$ and the structural tensors $F_\al$ are
covariantly constant with respect to the $(H,G)$-connection $D$ or
briefly $T$ and $F_\al$ are $D$-parallel, \ie $DT=0$ and
$DF_\al=0$. In the first subsection we give some interpretations
of these conditions.

\begin{prop}\label{th-DF-DT=0}
On an almost $(H,G)$-manifold the torsion $T$ of the
$(H,G)$-connection $D$ is $D$-parallel if and only if the
structural tensors $F_\al$ are $D$-parallel for all $\al=1,2,3$,
\ie $DT=0$ $\Leftrightarrow$ $DF_\al=0$, $\al=1,2,3$.
\end{prop}
\begin{pf}
Let $T$ be $D$-parallel. Then, from \eqref{QT} we obtain $DQ=0$.
After that, taking into account \eqref{FQ} and $DJ_\al=0$, we
obtain $DF_\al=0$, $\al=1,2,3$.

Vice versa, since we have $DF_\al=0$ and $DJ_\al=0$ for all
$\al=1,2,3$, from \eqref{Q} we obtain $DQ=0$. Hence $DT=0$ holds,
because of \eqref{T}. \hfill$\Box$\end{pf}

\begin{prop}\label{th-DF-Dt} On a locally conformal pseudo-hyper-K\"ahler manifold
the torsion $T$ of $D$ is $D$-parallel if and only if the 1-form
$\theta$ is $D$-parallel, \ie $DT=0$ $\Leftrightarrow$
$D\theta=0$.
\end{prop}
\begin{pf}
Let we have $DT=0$. Then the conditions $DF_\al=0$, $Dg=0$ and $D
J_\al=0$ for all $\al=1,2,3$ imply $D\theta=0$.

Vice versa, supposing conditions \eqref{F123} and $D\theta=0$
hold, we obtain $DF_\al=0$, $\al=1,2,3$, because of \eqref{DH}.
Then we have $DT=0$, according to \propref{th-DF-DT=0}. \phantom{}
\hfill$\Box$\end{pf}

If $\theta$ is $D$-parallel then we  obtain directly the following
\begin{cor}\label{cor-const}
On an almost $(H,G)$-manifold with $D$-parallel torsion of $D$ the
quantities $\theta(J_\al\Omega)$ are constant and, in particular,
$\theta(J_1\Omega)$ is zero.
\end{cor}

\subsection{Curvature properties of a $\W$-manifold with $D$-parallel torsion of $D$}
Here we characterize the locally conformal pseudo-hyper-K\"ahler
manifolds with $D$-parallel torsion of the $(H,G)$-connection $D$.
The main result is the following
\begin{thm}\label{thm-geom}
Every locally conformal pseudo-hyper-K\"ahler manifold $M$ with
$D$-parallel torsion of the complete natural connection $D$ is
$D$-flat and:
\begin{enumerate}\renewcommand{\labelenumi}{(\roman{enumi})}
    \item $M$ has a structure of a Lie group when is simply connected;
    \item $M$ is an isotropic pseudo-hyper-K\"ahler manifold;
    \item $M$ is scalar flat.
\end{enumerate}
\end{thm}

In the proof of the theorem we need to prove some lemmas.


\begin{lem}\label{prop-K-barV}
On a locally conformal pseudo-hyper-K\"ahler manifold with
$D$-par\-allel torsion of $D$ the curvature tensor
$K$ of $D$ has the form%
\be{KV} %
K=L-\frac{1}{256n^2}\{W+W\circ
J_1-W\circ J_2-W\circ J_3\}, %
\ee %
where %
\be{L}
L=R+\frac{1}{64n^2}\left(g\odot A\right),%
\ee
\be{A}
\begin{array}{ll}
& A=-\theta\otimes\theta+\theta\circ J_1\otimes\theta\circ
J_1+\theta\circ J_2\otimes\theta\circ J_2+\theta\circ
J_3\otimes\theta\circ
J_3\\[6pt]
&\phantom{
A=}-\theta(\Omega)g-\theta(J_2\Omega)g_2-\theta(J_3\Omega)g_3,
\end{array}
\ee
\be{barV}%
 W=g\odot C+\frac{1}{2}U, %
\ee%
\begin{equation}\label{C}
\begin{split}
& C= A+ B=2\theta\otimes\theta+2\sa\{\theta\circ
J_\al\otimes\theta\circ
J_\al\}\\[6pt]
&\phantom{ C= A+
B=}-\frac{5}{2}\theta(\Omega)g-\frac{3}{2}\theta(J_2\Omega)g_2
-\frac{3}{2}\theta(J_3\Omega)g_3.
\end{split}
\end{equation}
\end{lem}

\begin{pf}
Let the torsion $T$ of the $(H,G)$-connection $D$ be $D$-parallel
or equivalently $DF_\al=0$ ($\al=1,2,3$) hold. Then, for the
covariant derivative of
$F_\al$ ($\al=1,2,3$) with respect to  $\n$, we obtain %
\be{nFQ}%
\begin{split}
&\left(\n_x F_\al\right)(y,z,w)=F_\al(Q(x,y),z,w)\\[6pt]
&\phantom{\left(\n_x
F_\al\right)(y,z,w)=}+F_\al(y,Q(x,z),w)+F_\al(y,z,Q(x,w)).
\end{split}
\ee
Hence, having in mind \eqref{Ric-id-J}, \eqref{Q} and
\eqref{F123}, we obtain the following properties of the curvature
tensor $R$ for the Levi-Civita connection $\n$ of $g$:
\be{R} %
R-\ea R\circ J_\al=-\frac{1}{64n^2}\left\{ \left(g\odot
A\right)-\ea\left(g\odot A\right)\circ
J_\al\right\}, %
\ee %
where $A$ is determined by \eqref{A}.
It is obvious that $ A$ is symmetric, then $g\odot A$ is a
curvature-like tensor. Clearly, $R$ is a K\"ahler-like tensor if
and only if $ A=0$ or equivalently $\theta=0$. In other words, $R$
is a K\"ahler-like tensor if and only if $(M,H,G)$ is a
pseudo-hyper-K\"ahler manifold, which is known.

Now we transform \eqref{R} in the following form
\be{LL} %
R+\frac{1}{64n^2}\left(g\odot A\right)=\ea\Bigl\{R\circ
J_\al+\frac{1}{64n^2}\left(g\odot A\right)\circ
J_\al\Bigr\}, %
\ee %
which implies the property $L=\ea L\circ J_\al$ for the tensor $L$
from \eqref{L}. Therefore, $L$ is a K\"ahler-like tensor with
respect to each $J_\al$.

After that, applying \eqref{R} to  \eqref{K}, we
obtain%
\be{KL} K=L-\frac{1}{256n^2}\left\{\left(g\odot
C\right)+\sa\bigl\{\ea\left(g\odot C\right)\circ
J_\al\bigr\} +2U \right\}. \ee %

Obviously from \eqref{C} the tensor $ C$ is symmetric. Then
$g\odot C$ is a curvature-like tensor.

Let us remark that $W$, introduced by \eqref{barV}, has the K\"ahler-like property %
$W=\ea W\circ J_\al$.

Finally, from \eqref{KL} and \eqref{barV}, we obtain the form of
$K$ in terms of $L$ and $W$ given in \eqref{KV}.
\hfill$\Box$\end{pf}


\begin{lem}\label{prop-KLV=0}
On a locally conformal pseudo-hyper-K\"ahler manifold with
$D$-par\-allel torsion of $D$ the following tensors are zero:
\begin{enumerate}\renewcommand{\labelenumi}{(\roman{enumi})}
    \item the curvature tensor $K$ of the connection $D$ defined by \eqref{D};
    \item the tensor $W$ defined by \eqref{barV};
    \item the tensor $L$ defined by \eqref{L}.
\end{enumerate}
\end{lem}

\begin{pf}
It is known that $K$ is a K\"ahler-like tensor, \ie $K=\ea K\circ
J_\al$, because of $DJ_\al=0$ and \eqref{gJJ}. Moreover, we have
established that $L$, defined by \eqref{L}, is also a
K\"ahler-like tensor. Then the tensor $(W+W\circ J_1-W\circ
J_2-W\circ J_3)$ is also a K\"ahler-like tensor, according to
\eqref{KV}. Therefore, using  \propref{th-0}, we establish that
the tensors $K$, $L$ and $(W+W\circ J_1-W\circ J_2-W\circ J_3)$
are zero. It is easy to obtain that the last tensor in brackets is
zero if and only if $W=0$.

  \hfill$\Box$\end{pf}

\begin{pot}

Here we give some geometric interpretations of the results of
\lemref{prop-KLV=0}.

(i) The annulment of $K$ means that the manifold is $D$-flat.
Moreover, since $DT=0$ and $K=0$, according to the first Bianchi
identity of $K$ with torsion $T$, we obtain %
\be{sT} %
T\bigl(T(x,y),z\bigr)+T\bigl(T(y,z),x\bigr)+T\bigl(T(z,x),y\bigr)=0.%
\ee %
Thus, the expression $[\cdot,\cdot]=-T(\cdot,\cdot)$ is valid,
according to \cite{Hel}, \cite{KaTo}. Then \eqref{sT} becomes the
Jacobi identity and the manifold has a structure of a Lie group.

(ii) If we denote the trace
$\tau(L)=g^{is}g^{jk}L(e_i,e_j,e_k,e_s)$ of an arbitrary tensor
$L$ with respect to an arbitrary basis $\{e_i\}_{i=1}^{4n}$ then
we have for $W$ from \eqref{barV}
\[
\tau(W)=20n(4n-1)\theta(\Omega).
\]
Since $W=0$ according to \lemref{prop-KLV=0}, then
$\theta(\Omega)=0$ and because of \coref{cor-iK} we obtain that
the considered manifold is an isotropic pseudo-hyper-K\"ahler
manifold.

(iii) Since $L=0$ by \lemref{prop-KLV=0}, then according to
\eqref{L} the curvature tensor $R$,
the Ricci tensor $\rho$ and the scalar curvature $\tau$ have the form %
\be{R0} %
R=-\frac{1}{64n^2}\left(g\odot A\right), %
\ee
\[
\rho=\frac{1}{32n^2}\left\{(2n-1) A
-(2n+1)\theta(\Omega)g\right\},\quad
%
\tau=\frac{(2n+1)(1-4n)}{16n^2}\theta(\Omega). %
\]
%
Hence from $\theta(\Omega)=0$ it follows that %
\be{rhotau0} \tau=0, \qquad \rho=\frac{2n-1}{32n^2} A. \ee So,
the considered  manifold is scalar flat.

In that way, we complete the proof of \thmref{thm-geom}.
\hfill$\Box$\end{pot}

The following identity for the curvature tensor $R$ of the
Levi-Civita connection $\n$ on a complex manifold with Norden
metric $(M,J,g)$ is known from \cite{GrDj}:
\[
\begin{split}
&\mathop{\s}\limits_{x,y,z}
\bigl\{R(Jx,Jy,z,w)+R(x,y,Jz,Jw)\bigr. \\[6pt]
&\bigl.%
\phantom{\mathop{\s}\limits_{x,y,z}\bigl\{\bigr.}
+g\bigl(%
\left(\n_x J\right)y-\left(\n_y J\right)x,%
\left(\n_z J\right)w-\left(\n_w J\right)z%
\bigr)\bigr\}=0,
\end{split}
\]
where $\s$ denotes the cyclic sum by three arguments. Using the
last identity and the form \eqref{R0} of $R$, we obtain
\begin{cor}
On a locally conformal pseudo-hyper-K\"ahler manifold with
$D$-parallel torsion of $D$  it is valid the following identity
for $\al=2,3$
\begin{equation}\label{sgA}
    \mathop{\s}\limits_{x,y,z}\bigl(g_\al\odot
    A_\al\bigr)(x,y,z,w)=0,
\end{equation}
where  $A_\al=\theta\otimes\left(\theta\circ
J_\al\right)+\left(\theta\circ
J_\al\right)\otimes\theta$.$\hfill\Box$
\end{cor}

Now, having in mind \eqref{R0}, we obtain \be{nR0} \n_u
R=-\frac{1}{64n^2}\left(g\odot\n_u A\right) \ee and according to
\eqref{rhotau0} we have
 \be{nR0rho} \n_u
R=\frac{1}{2(1-2n)}\left(g\odot\n_u\rho\right). \ee
Since $\n_u\rho$ is a symmetric tensor, the last property implies
the second Bianchi identity without any additional condition. Then
a locally conformal pseudo-hyper-K\"ahler manifold $M$ with
$D$-parallel torsion  is Ricci-symmetric (\ie $\n\rho=0$) if and
only if $M$ is locally symmetric (\ie $\n R=0$). The same
proposition also follows  from \coref{cor-conf-flat}.


\section{An example of a 4-dimensional $\W$-manifold}

In \cite{GrMa} it is given an example of a 4-dimensional Lie group
as an almost $(H,G)$-manifold, which is an isotropic
pseudo-hyper-K\"ahler manifold. It belongs to the class $\W(J_1)$,
but is not a $\W$-manifold. The manifold is quasi-K\"ahlerian with
respect to $J_2$ and $J_3$ (\ie
$F_\al(x,y,z)+F_\al(y,z,x)+F_\al(z,x,y)=0$, $\al=2,3$).

In this section we construct an example of a 4-dimensional Lie
group as a $\W$-manifold.

Let $L$ be a 4-dimensional real connected Lie group, and
$\mathfrak{l}$ be its Lie algebra with a basis
$\{X_{1},X_{2},X_{3},X_{4}\}$.

Now we introduce an almost hypercomplex structure
$H=(J_1,J_2,J_3)$ by a standard way as in \eqref{Jdim4n} for
$k=0$:
\begin{equation}\label{Jdim4}
\begin{array}{llll}
J_1X_1=X_2,\quad & J_1X_2=-X_1,\quad & J_1X_3=-X_4,\quad
& J_1X_4=X_3,\\[4pt]
J_2X_{1}=X_{3}, \quad & J_2X_{2}=X_{4}, \quad & J_2X_{3}=-X_{1},
\quad &
J_2X_{4}=-X_{2}, \\[4pt]
J_3X_{1}=-X_{4},\quad & J_3X_{2}=X_{3}, \quad & J_3X_{3}=-X_{2},
\quad & J_3X_{4}=X_{1}.
\end{array}
\end{equation}

Let $g$ be a pseudo-Riemannian metric such that
\begin{equation}\label{g}
\begin{array}{c}
  g(X_1,X_1)=g(X_2,X_2)=-g(X_3,X_3)=-g(X_4,X_4)=1, \\[4pt]
  g(X_i,X_j)=0,\; i\neq j.
\end{array}
\end{equation}

Let us consider $(L,H,G)$ with the Lie algebra $\mathfrak{l}$
determined by the following nonzero commutators:
\begin{equation}\label{lie-w1-2}
\begin{array}{l}
\left
[X_{1},X_{4}\right]=[X_{2},X_{3}]=\lambda_{1}X_{1}+\lambda_{2}X_{2}+\lambda_{3}X_{3}+\lambda_{4}X_{4},\\[4pt]
\left[X_{1},X_{3}\right]=-[X_{2},X_{4}]=\lambda_{2}X_{1}-\lambda_{1}X_{2}+\lambda_{4}X_{3}-\lambda_{3}X_{4},
\end{array}
\end{equation}
where $\lambda_i\in\R$ ($i=1,2,3,4$).

We check that $J_\alpha$ ($\al=1,2,3$) are Abelian structures for
the Lie algebra $\mathfrak{l}$, \ie $[J_\alpha X_i,J_\alpha
X_j]=[X_i,X_j]$.

In \cite{Barb} five types of 4-dimensional Lie algebras which
admit integrable invariant hypercomplex structures are classified.
Obviously, the introduced Lie algebra $\mathfrak{l}$ in
\eqref{lie-w1-2} with the hypercomplex structure from
\eqref{Jdim4} is of the third type in \cite{Barb} when
$\lambda_i\neq 0$ for some $i=1,2,3,4$. Namely, the Lie algebra
$\mathfrak{l}$ is isomorphic to $\mathfrak{aff}(\mathbb{C})$.

Then we prove the following
\begin{thm}\label{HG}
Let $(L,H,G)$ be the almost $(H,G)$-manifold, determined by
\eqref{Jdim4}, \eqref{g} and \eqref{lie-w1-2}. Then:
\begin{enumerate}\renewcommand{\labelenumi}{(\roman{enumi})}
    \item it belongs to the class of the locally conformal
        pseudo-hyper-K\"ahler manifolds;
    \item  it is locally conformally flat and the
        curvature tensor $R$ has the form
        \begin{equation}\label{Rform3}
            \begin{array}{l}
                R=-\frac{1}{2}(g\odot\rho)+\frac{\tau}{12}(g\odot
                g);
            \end{array}
        \end{equation}
    \item it is $D$-flat with non-$D$-parallel torsion of $D$.
\end{enumerate}
\end{thm}

\begin{pf}
According to (\ref{g}), (\ref{lie-w1-2}) and the well-known
property of a Levi-Civita connection $\n$,
we obtain the nonzero components of $\nabla$ generated by $g$, as
follows:

\begin{equation}\label{nabla3}
\begin{array}{ll}
\nabla_{X_{1}}X_{1} = \nabla_{X_{2}}X_{2} = \lambda_{2}X_{3} +
\lambda_{1}X_{4},
\\[4pt]
\nabla_{X_{1}}X_{3} = \nabla_{X_{4}}X_{2} = \lambda_{2}X_{1} -
\lambda_{3}X_{4},
\\[4pt]
\nabla_{X_{1}}X_{4} = -\nabla_{X_{3}}X_{2} = \lambda_{1}X_{1} +
\lambda_{3}X_{3},
\\[4pt]
\nabla_{X_{2}}X_{3} =- \nabla_{X_{4}}X_{1} =\lambda_{2}X_{2} +
\lambda_{4}X_{4},
\\[4pt]
\nabla_{X_{2}}X_{4} = \nabla_{X_{3}}X_{1} = \lambda_{1}X_{2} -
\lambda_{4}X_{3},
\\[4pt]
\nabla_{X_{3}}X_{3} = \nabla_{X_{4}}X_{4} = -\lambda_{4}X_{1} -
\lambda_{3}X_{2}.
\end{array}
\end{equation}

The components of $\nabla J_\alpha$ follow from the last
equalities and \eqref{Jdim4}. Then, having in mind (\ref{g}), we
obtain the following nonzero components
$(F_\alpha)_{ijk}=F_\alpha(X_{i},X_{j},X_{k})=g\left(\left(\nabla_{X_i}J_\alpha\right)X_j,X_k\right)$
of the tensors $F_\alpha$, $\al=1,2,3$:
\begin{subequations}\label{lambdi}
\begin{equation}
\begin{array}{rl}
&\lambda_{1}=(F_1)_{113}=(F_1)_{124}=-(F_1)_{131}=-(F_1)_{142}\\[4pt]
&\phantom{\lambda_{1}}=-(F_1)_{214}=(F_1)_{223}=-(F_1)_{232}=(F_1)_{241}\\[4pt]
&\phantom{\lambda_{1}}=(F_2)_{112}=(F_2)_{121}=(F_2)_{134}=(F_2)_{143}=\frac{1}{2}(F_2)_{222}\\[4pt]
&\phantom{\lambda_{1}}=\frac{1}{2}(F_2)_{244}=(F_2)_{314}
=-(F_2)_{323}=-(F_2)_{332}=(F_2)_{341}\\[4pt]
&\phantom{\lambda_{1}}=-\frac{1}{2}(F_3)_{111}=-\frac{1}{2}(F_3)_{144}=-(F_3)_{212}=-(F_3)_{221}=(F_3)_{234}\\[4pt]
&\phantom{\lambda_{1}}=(F_3)_{243}=(F_3)_{313}
=(F_3)_{324}=(F_3)_{331}=(F_3)_{342},\\[4pt]
\end{array}
\end{equation}
\begin{equation}
\begin{array}{rl}
&\lambda_{2}=-(F_1)_{114}=(F_1)_{123}=-(F_1)_{132}=(F_1)_{141}\\[4pt]
&\phantom{\lambda_{2}}=-(F_1)_{213}=-(F_1)_{224}=(F_1)_{231}=(F_1)_{242}\\[4pt]
&\phantom{\lambda_{2}}=\frac{1}{2}(F_2)_{111}=\frac{1}{2}(F_2)_{133}=(F_2)_{212}=(F_2)_{221}=(F_2)_{234}\\[4pt]
&\phantom{\lambda_{2}}=(F_2)_{243}=-(F_2)_{414}=(F_2)_{423}=(F_2)_{432}
=-(F_2)_{441}\\[4pt]
&\phantom{\lambda_{2}}=(F_3)_{112}=(F_3)_{121}=-(F_3)_{134}=-(F_3)_{143}=\frac{1}{2}(F_3)_{222}\\[4pt]
&\phantom{\lambda_{2}}=\frac{1}{2}(F_3)_{233}=-(F_3)_{413}=-(F_3)_{424}=-(F_3)_{431}
=-(F_3)_{442},\\[4pt]
\end{array}
\end{equation}
\begin{equation}
\begin{array}{rl}
&\lambda_{3}=(F_1)_{313}=(F_1)_{324}=-(F_1)_{331}=-(F_1)_{342}\\[4pt]
&\phantom{\lambda_{3}}=(F_1)_{414}=-(F_1)_{423}=(F_1)_{432}=-(F_1)_{441}\\[4pt]
&\phantom{\lambda_{3}}=(F_2)_{114}=-(F_2)_{123}=-(F_2)_{132}=(F_2)_{141}=-(F_2)_{312}\\[4pt]
&\phantom{\lambda_{3}}=-(F_2)_{321}=-(F_2)_{334}=-(F_2)_{343}=-\frac{1}{2}(F_2)_{422}
=-\frac{1}{2}(F_2)_{444}\\[4pt]
&\phantom{\lambda_{3}}=(F_3)_{113}=(F_3)_{124}=(F_3)_{131}=(F_3)_{142}=-\frac{1}{2}(F_3)_{322}\\[4pt]
&\phantom{\lambda_{3}}=-\frac{1}{2}(F_3)_{333}=(F_3)_{412}=(F_3)_{421}=-(F_3)_{434}
=-(F_3)_{443},\\[4pt]
\end{array}
\end{equation}
\begin{equation}
\begin{array}{rl}
&\lambda_{4}=(F_1)_{314}=-(F_1)_{323}=(F_1)_{332}=-(F_1)_{341}\\[4pt]
&\phantom{\lambda_{4}}=-(F_1)_{413}=-(F_1)_{424}=(F_1)_{431}=(F_1)_{442}\\[4pt]
&\phantom{\lambda_{4}}=-(F_2)_{214}=(F_2)_{223}=(F_2)_{232}=-(F_2)_{241}=-\frac{1}{2}(F_2)_{311}\\[4pt]
&\phantom{\lambda_{4}}=-\frac{1}{2}(F_2)_{333}=-(F_2)_{412}=-(F_2)_{421}=-(F_2)_{434}
=-(F_2)_{443}\\[4pt]
&\phantom{\lambda_{4}}=-(F_3)_{213}=-(F_3)_{224}=-(F_3)_{231}=-(F_3)_{242}=-(F_3)_{312}\\[4pt]
&\phantom{\lambda_{4}}=-(F_3)_{321}=(F_3)_{334}=(F_3)_{343}=\frac{1}{2}(F_3)_{411}
=\frac{1}{2}(F_3)_{444}.\\[4pt]
\end{array}
\end{equation}
\end{subequations}

Using \eqref{lambdi}, we establish that \eqref{cl-H} and
\eqref{cl-N} are satisfied. Therefore the mani\-fold $(L,H,G)$
belongs to the class $\W$.

After that, by \eqref{theta-al}, \eqref{titi-all} and
\eqref{Jdim4}, we obtain the components
$\left(\theta\right)_k=\theta(X_k)$ of the 1-form $\theta$:
\begin{equation}\label{theta123}
\begin{array}{llll}
\left(\theta\right)_{1}=4\lambda_{4},  \quad &
\left(\theta\right)_{2}=4\lambda_{3}, \quad &
\left(\theta\right)_{3}=-4\lambda_{2}, \quad &
\left(\theta\right)_{4}=-4\lambda_{1}.
\end{array}
\end{equation}

(i) We verify, using \eqref{lie-w1-2} and \eqref{theta123}, that
$\theta$ is closed (\ie $\D\theta=0$) and therefore $(L,H,G)$ is a
locally conformal pseudo-hyper-K\"ahler manifold.

(ii) It is clear, because of (i) and \coref{cor-conf-flat}, that
$(L,H,G)$ is locally conformally flat. Thus, the corresponding
Weyl tensor vanishes and conse\-quent\-ly the curvature tensor $R$
has the form in \eqref{Rform3}.

(iii)  Using \eqref{nabla3}, \eqref{Jdim4} and the components of
$\nabla J_\alpha$, we obtain the components of the
$(H,G)$-connection $D$, determined by \eqref{D}:
\begin{equation}
\begin{array}{ll}
D_{X_1}X_1=-\frac{1}{2}\lambda_3X_2,\qquad &
D_{X_1}X_2=\frac{1}{2}\lambda_3X_1,\\[4pt]
D_{X_1}X_3=-\frac{1}{2}\lambda_3X_4,\qquad &
D_{X_1}X_4=\frac{1}{2}\lambda_3X_3,\\[4pt]
D_{X_2}X_1=\frac{1}{2}\lambda_4X_2,\qquad &
D_{X_2}X_2=-\frac{1}{2}\lambda_4X_1,\\[4pt]
D_{X_2}X_3=\frac{1}{2}\lambda_4X_4,\qquad &
D_{X_2}X_4=-\frac{1}{2}\lambda_4X_3,\\[4pt]
D_{X_3}X_1=\frac{1}{2}\lambda_1X_2,\qquad &
D_{X_3}X_2=-\frac{1}{2}\lambda_1X_1,\\[4pt]
D_{X_3}X_3=\frac{1}{2}\lambda_1X_4,\qquad &
D_{X_3}X_4=-\frac{1}{2}\lambda_1X_3,\\[4pt]
D_{X_4}X_1=-\frac{1}{2}\lambda_2X_2,\qquad &
D_{X_4}X_2=\frac{1}{2}\lambda_2X_1,\\[4pt]
D_{X_4}X_3=-\frac{1}{2}\lambda_2X_4,\qquad &
D_{X_4}X_4=\frac{1}{2}\lambda_2X_3.\\[4pt]
\end{array}
\end{equation}
Hence we establish that the corresponding curvature tensor $K$
vanishes, therefore $(L,H,G)$ is $D$-flat.

We verify immediately that $(L,G,H)$ has no $D$-parallel torsion
of $D$. If we put the condition $DT=0$, then the Lie algebra
$\mathfrak{l}$ from \eqref{lie-w1-2} becomes Abelian. $\phantom{}
\hfill \Box$
\end{pf}

According to (\ref{g}) and (\ref{nabla3}), the curvature tensor
$R$ has the following nonzero components
$R_{ijkl}=R(X_{i},X_{j},X_{k},X_{l})$:
\begin{equation} \label{R3}
\begin{array}{c}
\begin{array}{lll}
R_{1221} = \lambda_{1}^{2} + \lambda_{2}^{2}, \qquad & %
R_{1331} = \lambda_{4}^{2} - \lambda_{2}^{2}, \qquad & %
R_{1441} = \lambda_{4}^{2} - \lambda_{1}^{2},
\\[4pt]
R_{2332} = \lambda_{3}^{2} - \lambda_{2}^{2},\qquad & %
R_{2442} = \lambda_{3}^{2} - \lambda_{1}^{2}, \qquad & %
R_{3443} = -\lambda_{3}^{2} - \lambda_{4}^{2},\\[4pt]
\end{array}
\\[4pt]
\begin{array}{ll}
R_{1341}=R_{2342} = -\lambda_{1}\lambda_{2}, \qquad & %
R_{2132}=-R_{4134} = -\lambda_{1}\lambda_{3},
\\[4pt]
R_{1231}=-R_{4234} = \lambda_{1}\lambda_{4}, \qquad & %
R_{2142}=-R_{3143} = \lambda_{2}\lambda_{3},
\\[4pt]
R_{1241}=-R_{3243} = -\lambda_{2}\lambda_{4}, \qquad & %
R_{3123}=R_{4124} = \lambda_{3}\lambda_{4},
\end{array}
\end{array}
\end{equation}
and the rest are determined by \eqref{R3} and properties
\eqref{curv}.

Hence, (\ref{R3}) implies the components of the Ricci tensor
$\rho$ and the value of the scalar curvature $\tau$:
\begin{equation} \label{Ricci3}
\begin{array}{ll}
\rho_{11}=2\big( \lambda_{1}^{2} + \lambda_{2}^{2} - \lambda_{4}^{2} \big), \; & %
\rho_{12}=\rho_{21}=-2\lambda_{3}\lambda_{4}, \;\; %
\rho_{23}=\rho_{32}=2\lambda_{1}\lambda_{4},
\\[4pt]
\rho_{22}=2\big( \lambda_{1}^{2} + \lambda_{2}^{2} - \lambda_{3}^{2} \big), \; & %
\rho_{13}=\rho_{31}=-2\lambda_{1}\lambda_{3},\;\; %
\rho_{24}=\rho_{42}=-2\lambda_{2}\lambda_{4},
\\[4pt]
\rho_{33}=2\big( \lambda_{4}^{2} + \lambda_{3}^{2} - \lambda_{2}^{2} \big),\; &%
\rho_{14}=\rho_{41}=2\lambda_{2}\lambda_{3},\;\; \phantom{-}%
\rho_{34}=\rho_{43}=-2\lambda_{1}\lambda_{2},
\\[4pt]
\rho_{44}=2\big( \lambda_{4}^{2} + \lambda_{3}^{2} - \lambda_{1}^{2} \big),\;\; &%
\phantom{\rho_{44}..\rho_{44}}\tau=6\big(\lambda_{1}^{2}+\lambda_{2}^{2}-\lambda_{3}^{2}-\lambda_{4}^{2}\big).
\end{array}
\end{equation}

\begin{cor}\label{th-iK}
The manifold $(L,H,G)$ is scalar flat and isotropic
hyper-K\"ahlerian  if and only if the structural constants from
\eqref{lie-w1-2} satisfy the following condition
\begin{equation}\label{llll}
\begin{array}{l}
\lambda_{1}^{2}+\lambda_{2}^{2}-\lambda_{3}^{2}-\lambda_{4}^{2}=0.
\end{array}
\end{equation}
\end{cor}
\begin{pf}
It follows immediately from the square norms of $\nabla J_\alpha$: %
\[
-2\nJ{1}=\nJ2=\nJ{3}=16\left(\lambda_1^2+\lambda_2^2-\lambda_3^2-\lambda_4^2\right)
\]
and the last equation of \eqref{Ricci3}. \hfill$\Box$\end{pf}

\section*{Acknowledgements} The author would like to express his deep gratitude to
his doctoral advisor Kostadin Gribachev for the latter's
invaluable help, patience and support over the last 22 years and
in  connection with his 70th anniversary.

The author wishes to thank both Stefan Ivanov for several useful
discussions about this work and the referee, who contributed to
improve considerably the final presentation of the manuscript.

\end{document}